\documentstyle[12pt]{article}
\begin{document}

\baselineskip 20pt

\vspace{.5in}

\begin{center}

{\Large \bf
A special class of solutions\\ of \\
the truncated Hill's equation\\
}
\vskip .8in
{\bf
Shang Yuan Ren\\
Department of Physics, Peking University \\
Beijing 100871, People's Republic of China\\
}
\vbox{}
\end{center}
\newpage
\begin{center}
{\bf \large Abstract}
\end{center}
\par
This work investigates the existence and properties of a certain class of
solutions of the Hill's equation truncated in the interval $[\tau, \tau + L]$
- where $L = Na$, $a$ is the period of the coefficients in Hill's equation,
$N$ is a positive integer and $\tau$ is a real number. It is found that
the truncated Hill's equation has two different types of solutions
which vanish at the truncation boundaries $\tau$ and $\tau + L$:
There are always $N-1$ solutions in each stability interval of
Hill's equation, whose eigen value
is dependent on the truncation length $L$ but not on the truncation boundary
$\tau$;
There is always one and only one solution in each finite conditional
instability interval of Hill's equation, whose eigen value might be
dependent on the boundary location $\tau$ but not on the truncation length $L$.
\par
The results obtained are applied to the physics problem on the electronic
states in one dimensional crystals of finite length. It significantly
improves many known results and also provides more new understandings
on the physics in low dimensional systems.
\vskip .8in
\newpage
\section{Introduction}
Hill's equation is a class of second-order ordinary differential equations
with coefficients of period $a$.
The properties of solutions of Hill's equation has been
investigated in detail by many mathematicians\cite{mag,eas}.
It is well known that the eigenvalues and eigenfunctions of an
ordinary second order differential equation depend on both the region in which
the equation is given and the boundary conditions.
If the equation is given in a finite interval and we are interested in
a special class of the solutions - they vanish
at the boundaries - we expect that both the eigenvalues and the
eigenfunctions are dependent on the boundary locations and the properties of
the finite interval.
For an ordinary differential equation with periodic coefficients, in
general it is usually
expected that both the eigenvalues and the eigenfunctions are dependent on
{\it both the boundary locations and the interval length}.\par
By using the general theory of Hill's equation, this work investigates
the existence and properties of this class of solutions of
Hill's equation truncated in the interval $[\tau, \tau + L]$ - where
$L = Na$, $a$ is the period of the coefficients in Hill's equation and
$N$ is a positive integer and $\tau$ is a real number. It is found that
the truncated Hill's equation have two different types of solutions
which vanish at the truncation boundaries $\tau$ and $\tau + L$:
There are always $N-1$ solutions in each stability interval of
Hill's equation, whose eigen value
is dependent on the truncation length $L$ but not on the truncation boundary
$\tau$;
There is always one and only one solution in each finite conditional
instability interval, whose eigen value might be dependent on the
boundary location $\tau$ but not on the truncation length $L$.
\par
The Schr$\ddot{\rm o}$dinger differential equation on
the electronic states in one dimensional crystals\cite{kit}
can be considered as a special case of Hill's equation.
Due to the fact that any real crystal always has a finite size,
the electronic states in a real crystal of finite size can be better described
as solutions of a special case of Hill's equation truncated in a finite length
which vanish at the boundaries.
The general results obtained here are applied to
the physics problem on the electronic
states in one dimensional crystals of finite length. It significantly
improves many known results and also provides many new understandings
on the physics in low dimensional systems.\par
\section{Hill's equation}
We consider the Hill's equation having the following form:
\begin{equation}
\{ p(x)y'(x) \}' +  \{ \lambda  s(x) - q(x) \} y(x) = 0,
~~~~~~~~~-\infty~<~x~<~+\infty
\end{equation}
where $\lambda$ is a real parameter,
 $p(x), ~q(x)$, and $s(x)$ are real-valued and have the same period $a$.
It is assumed that $p(x)$ is continuous and nowhere zero and that
$p'(x)$, $s(x)$ and $q(x)$ are piecewise
 continuous with period $a$ and there
is a constant $s > 0$ such that $s(x) \geq s$\cite{eas}.
\par
A periodic eigenvalue problem and a semi-periodic
eigen value problem play a significant role in the theory of
Hill's equation (1)\cite{eas}:\\
(i). The periodic eigenvalue problem comprises (1), considered to hold
in $[0, a]$ and the periodic conditions:
$$
y(a) = y(0),~~~y'(a) = y'(0).
$$
The eigenfunctions are denoted by $\zeta_n(x)$ and eigenvalues by $\nu_n$
(n = 0, 1,....).
\\
(ii). The semi-periodic eigenvalue problem comprises (1), considered to hold
in $[0, a]$ and the semi-periodic conditions:
$$
y(a) = -y(0),~~~y'(a) = -y'(0).
$$
The eigenfunctions are denoted by $\xi_n(x)$ and eigenvalues by $\mu_n$
(n = 0, 1,....).
\par
The numbers $\nu_n$ and $\mu_n$ occur in the order\cite{eas}
$$
\nu_0 < \mu_0 \leq \mu_1 <
\nu_1 \leq
\nu_2 < \mu_2 \leq \mu_3 < \nu_3  \leq \nu_4 < .....
$$
\par
Let $\phi_1(x, \lambda)$ and $\phi_2(x,\lambda)$ be the linearly independent
solutions of
(1), such that
$
\phi_1(0, \lambda) = 1,~
\phi_1'(0, \lambda) = 0;~~
\phi_2(0, \lambda) = 0,~
\phi_2'(0, \lambda) = 1,
$
the real number $D(\lambda)$ defined by
$
D (\lambda) = \phi_1(a, \lambda) + \phi_2'(a, \lambda)
$
is called the $discriminant$ of (1) and it plays a significant role in
determining the properties of the solutions of (1). \par
$D(\lambda)$ is an analytical function of $\lambda$\cite{eas}.
The values
of $\lambda$ for which $|D(\lambda)| < 2$ form an open set on the real
$\lambda$-axis. This set can be expressed as the union of countable collection
of disjoint open intervals
$(\nu_{2m},\mu_{2m})$ and $(\mu_{2m+1},\nu_{2m+1})$.
Equation (1) is stable when $\lambda$ lies in these intervals thus the
intervals are called the stability intervals of (1). The set on the real
$\lambda$-axis with the stability intervals excluded can be expressed as the
union of countable collection
of disjoint intervals in which $ |D(\lambda)| \geq 2$:
$(-\infty, \nu_0]$, $[\mu_{2m},\mu_{2m+1}]$ and $[\nu_{2m+1},\nu_{2m+2}]$.
Those intervals are called conditional instability intervals.
Thus as $\lambda$ increases from $- \infty$ to $+ \infty$, the conditional
instability intervals and the stability intervals occur alternatively.
\section{The truncated Hill's equation}
\par
A special class of solutions of the truncated Hill's equation are the main
interest in this paper.
We are trying to find out the eigenvalues $\Lambda$ and eigenfunctions
$y(x)$, which are solutions of
\begin{equation}
\{ p(x)y'(x) \}' +  \{ \Lambda  s(x) - q(x) \} y(x) = 0,
~~~~~~~~~\tau~\leq~x~\leq~\tau + L
\end{equation}
under the following boundary condition:
\begin{equation}
y( \tau, \Lambda) = y( \tau +L, \Lambda) = 0,
\end{equation}
where $\tau$ is a real number, $ L = Na$ and $N$ is a positive integer.
\par
Suppose $y_1(x, \lambda)$ and $y_2(x, \lambda)$ are two linearly independent
solutions of (1), in general, the solution of (2) which satisfies
(3) if it exists, can be expressed as
\begin{equation}
y( x, \Lambda) =  c_1 y_1(x, \Lambda) + c_2 y_2(x, \Lambda).
~~~~~~~ ~~\tau~\leq~x~\leq~\tau +L
\end{equation}
In the following we try to find out the eigenvalues of (2) under the
boundary condition (3) by assuming all non-trivial solutions of (1) are known.
After each eigenvalue is found, the corresponding eigenfunction can be
obtained easily.\par
In principle, we should consider
solutions of Equations (2) and (3) for $\lambda$ in $(-\infty, +\infty)$.
However, according to the Theorem 3.2.2 of \cite{eas}, there is not
a nontrivial solution of (2) satisfying (3) for
$\Lambda$ in $(- \infty, \nu_0]$.
Thus we need only to consider $\lambda$ in $(\nu_0, +\infty)$.
\par
For $\lambda$ in $(\nu_0, +\infty)$, $D(\lambda) $ may have five different
cases depending on the value of $\lambda$.
Correspondingly, the two
linearly independent solutions of (1) are also different:
\\
Case A. $|D(\lambda)| < 2$:
\par
This is in the stability intervals of (1).
Two linearly independent solutions of (1) can be expressed as\cite{eas}
$$
y_1(x, \lambda) = e^{i \alpha(\lambda) x} p_1(x, \lambda),~~
y_2(x, \lambda) = e^{- i \alpha(\lambda) x} p_2(x, \lambda),
$$
where $\alpha(\lambda)$ is a real number depending on $\lambda$ and
\begin{eqnarray*}
0 < \alpha(\lambda) a < \pi
\end{eqnarray*}
and $p_1(x, \lambda)$ and $p_2(x, \lambda)$ have period $a$.
All $\alpha(\lambda)$ and $p_1(x, \lambda)$ and $p_2(x, \lambda)$ are functions
of $\lambda$.
\par
From (3) and (4) we have
$$
c_1 e^{i \alpha \tau} p_1(\tau, \lambda) +
c_2 e^{- i \alpha \tau} p_2( \tau, \lambda)
 = 0,
$$
and
$$
c_1 e^{i \alpha (\tau + L)} p_1(\tau + L, \lambda) +
c_2 e^{-i \alpha (\tau + L) } p_2(\tau +L, \lambda) = 0.
$$
Due to that
$$
p_1(\tau, \lambda) = p_1(\tau + L, \lambda), ~~~
p_2(\tau, \lambda) = p_2(\tau + L, \lambda),
$$
we obtain either
$$
e^{i \alpha(\lambda) L} - e^{- i \alpha(\lambda) L} = 0 .
\eqno(A.1)
$$
or
$$
c_1 p_1(\tau, \lambda) = 0~~~ and ~~~~ c_2 p_2( \tau, \lambda) = 0.
\eqno(A.2)
$$
We can easily prove that neither $p_1(\tau, \lambda)$ nor $p_2(\tau, \lambda)$
can be zero. Suppose $p_1(\tau, \lambda) = 0$, we have
$p_1(\tau +a, \lambda) = 0$
and thus $y_1(\tau, \lambda) = y_1(\tau +a, \lambda) = 0$. Then
according to the theorem 3.1.3 of \cite{eas}, $\lambda$ must be in the range of
$[\mu_{2m}, \mu_{2m+1}]$
or $[\nu_{2m+1}, \nu_{2m+2}]$, in which $ | D(\lambda) | \geq 2$. This is in
contradictory to the condition $ | D(\lambda)| < 2 $ here. Similarly
$p_2(\tau, \lambda)$ can not be zero. Thus we must have $c_1 = c_2 = 0$
and the non-trivial solutions obtained from (A.2) do not exist.\par
Thus the existence of non-trivial solution of (2) and (3) requires
$$
\alpha(\lambda) L = j \pi,~~~~~~~~~~~~ j~=~1,~2,~3,......N-1.
$$
Therefore for each stability interval, there are $N-1$ values of $\Lambda_j$, where
$j = 1, 2, ....N-1$, for which
\begin{equation}
\alpha(\Lambda_j) = j~ \pi/L .
\end{equation}
Each eigen value for this case, is a function of $L$, the truncation length.
But they all do not depend on the location of the truncation boundary $\tau$.
\\
Case B. $D(\lambda) = 2$.
\par
There could be two possibilities\cite{eas}:
\\
B1. $\phi_2(a, \lambda)$ and $ \phi_1'(a, \lambda)$ are both zero:
$$
 \phi_2(a, \lambda)  ~=~\phi_1'(a, \lambda)~=~0.
$$
This corresponds to $\nu_{2m+1} = \nu_{2m+2}$.
Two linearly independent solutions of (1) can be expressed as\cite{eas}
$$
y_1(x, \lambda) = p_1(x, \lambda),~~
y_2(x, \lambda) = p_2(x, \lambda),
$$
and $p_1(x, \lambda)$ and $p_2(x, \lambda)$ have period $a$.
It is always possible to choose $c_1$ and $c_2$, which are not both zero, to
make $y(\tau, \lambda) = c_1 y_1(\tau, \lambda) + c_2 y_2(\tau, \lambda) = 0$
and naturally we have $y(\tau + L, \lambda) = 0$.
Thus we have a solution $ \Lambda~=~\lambda $  of (2) and (3) in this case
and the corresponding eigenfunction is a periodic function.
\par
Obviously for this case, we always have
$$
y(\tau, \lambda)~=~y(\tau + a, \lambda) ~=~ 0.
$$
B2. $\phi_2(a, \lambda)$ and $ \phi_1'(a, \lambda)$ are not both zero.
\par
Two linearly independent solutions of (1) can be expressed as\cite{eas}
$$
y_1(x, \lambda) = p_1(x, \lambda),~~
y_2(x, \lambda) = x p_1(x, \lambda) + p_2(x, \lambda),
$$
and $p_1(x, \lambda)$ and $p_2(x, \lambda)$ have period $a$.
\par
Equations (3) and (4) become
$$
c_1 p_1( \tau, \lambda) + c_2 ( \tau p_1 ( \tau, \lambda) +
p_2 ( \tau, \lambda)) = 0,
$$
and
$$
c_1 p_1( \tau + L, \lambda) +
c_2 ( (\tau + L) p_1 ( \tau + L, \lambda) + p_2 ( \tau + L, \lambda)) = 0.
$$
Thus, we must have,
$$
c_2 L p_1 ( \tau, \lambda) = 0,
$$
and
$$
c_1 p_1( \tau, \lambda) +
c_2 ( \tau p_1 ( \tau, \lambda) + p_2 ( \tau, \lambda)) = 0.
$$
These two lead to
$$
p_1( \tau, \lambda) = 0, ~~~~~~~~~~~and ~~~~~~~~~~~ c_2 = 0.
\eqno (B2.1)
$$
Thus if we have a solution for this case, we must have (B2.1).  \\
Case C. $ D(\lambda) > 2$.
\par
Two linearly independent solutions of (1) can be expressed as\cite{eas}
$$
y_1(x, \lambda) = e^{\beta(\lambda) x} p_1(x, \lambda),~~
y_2(x, \lambda) = e^{-\beta(\lambda) x} p_2(x, \lambda),
$$
where $\beta(\lambda)$ is a real non-zero number depending on $\lambda$ and
$p_1(x, \lambda)$ and $p_2(x, \lambda)$ have period $a$.
\par
Equations (3) and (4) lead to
$$
c_1 e^{ \beta \tau} p_1(\tau, \lambda) +
c_2 e^{- \beta \tau} p_2( \tau, \lambda) = 0,
$$
and
$$
c_1 e^{ \beta (\tau + L)} p_1(\tau + L, \lambda) +
c_2 e^{- \beta (\tau + L) } p_2( \tau + L, \lambda) = 0.
$$
Due to
$$
p_1(\tau,\lambda) = p_1(\tau + L,\lambda), ~~~
p_2(\tau,\lambda) = p_2(\tau + L,\lambda), ~~~
$$
we must have
$$
c_1 p_1(\tau, \lambda) = 0, ~~ and ~~ c_2 p_2( \tau, \lambda) = 0.
$$
By the Sturm Separation Theorem\cite{eas2}, the zeros of
$p_1(x, \Lambda) $ are separated from the zeros of $ p_2(x, \Lambda) $.
Thus $p_1(\tau, \lambda) $ and $ p_2( \tau, \lambda) $ can not be zero
simultaneously, either $c_1$ or $c_2$ must be zero.
Whether we have a solution or not in this case is dependent on
whether either one of
$$
p_1(\tau, \lambda) = 0 ~~~~~~~and~~~~~~~~c_2 = 0
\eqno(C.1)
$$
and
$$
p_2( \tau, \lambda) = 0 ~~~~~~~and~~~~~~~~c_1 = 0
\eqno(C.2)
$$
is true.
From the discussion on (B1), (B2), (C), we can see that in the
conditional instability intervals in which $D(\lambda) \geq 2$, if we have
a solution of (2) and (3), we must have either $\nu_{2m+1} = \nu_{2m+2}$
or one of (B2.1), (C.1) and (C.2) when $\nu_{2m+1} < \nu_{2m+2}$.
Thus we always have $y(\tau~+~a, \Lambda) = 0$ if
$y(\tau, \Lambda) =0$.
Therefore for conditional instability intervals $D(\lambda) \geq 2$,
the following equation is a necessary
condition for having a solution of Equations (2) and (3):
\begin{equation}
y(\tau~+~a, \Lambda)~ =~y(\tau, \Lambda)~=~0.
\end{equation}
Case D. $ D(\lambda) = -2$.\\
There could be two possibilities\cite{eas}:
\par
D1. $\phi_2(a, \lambda)$  and $\phi_1'(a, \lambda)$ are both zero:
$$
\phi_2(a, \lambda) = \phi_1'(a, \lambda) = 0.
$$
This corresponds to $\mu_{2m} = \mu_{2m+1}$.
Two linearly independent solutions of (1) can be expressed as\cite{eas}
$$
y_1(x, \lambda) = s_1(x, \lambda),~~
y_2(x, \lambda) = s_2(x, \lambda),
$$
and $s_1(x, \lambda)$ and $s_2(x, \lambda)$ have semi-period $a$:
$s_i(x+a) = - s_i(x)$.
It is always possible to choose $c_1$ and $c_2$, they are not both zero, to
make $y(\tau, \lambda) = c_1 y_1(\tau, \lambda) + c_2 y_2(\tau, \lambda) = 0$
and we naturally have $y(\tau, \lambda) = y(\tau + L, \lambda) = 0$.
Thus we always have a solution of (2) and (3) in this case.
The corresponding eigenfunction is a semiperiodic function.
\par
Obviously for this case, we have
$$
y(\tau, \lambda)~=~y(\tau + a, \lambda) ~=~ 0.
$$
D2. $\phi_2(a, \lambda)$ and $ \phi_1'(a, \lambda)$ are not both zero.\par
Two linearly independent solutions of (1) can be expressed as\cite{eas}
$$
y_1(x, \lambda) = s_1(x, \lambda),~~
y_2(x, \lambda) = x~ s_1(x, \lambda) + s_2(x, \lambda),
$$
and $s_1(x, \lambda)$ and $s_2(x, \lambda)$ have semi-period $a$.
\par
Equations (3) and (4) become
$$
c_1 s_1( \tau, \lambda) +
c_2 ( \tau s_1 ( \tau, \lambda) + s_2 ( \tau, \lambda)) = 0,
$$
and
$$
c_1 s_1( \tau + L, \lambda) +
c_2 ( (\tau + L) s_1 ( \tau + L, \lambda) + s_2 ( \tau + L, \lambda)) = 0.
$$
Thus, we must have,
$$
c_2 L s_1 ( \tau, \lambda) = 0,
$$
and
$$
c_1 s_1( \tau, \lambda) +
c_2 ( \tau s_1 ( \tau, \lambda) + s_2 ( \tau, \lambda)) = 0.
$$
These two lead to
$$
s_1( \tau, \lambda) = 0~~~~~~~~~ and~~~~~~~~~~~c_2 = 0.
\eqno (D2.1)
$$
Case E. $ D(\lambda) < - 2$:\\
Two linearly independent solutions of (1) can be expressed as\cite{eas}
$$
y_1(x, \lambda) = e^{\beta(\lambda)  x} s_1(x, \lambda),~~
y_2(x, \lambda) = e^{-\beta(\lambda) x} s_2(x, \lambda),
$$
where $\beta(\lambda)$ is a non-zero real number depending on $\lambda$ and
$s_1(x, \lambda)$ and $s_2(x, \lambda)$ are semiperiodic functions.
\par
Equations (3) and (4) lead to
$$
c_1 e^{ \beta(\lambda) \tau} s_1(\tau, \lambda) +
c_2 e^{-\beta(\lambda) \tau} s_2( \tau, \lambda) = 0,
$$
and
$$
c_1 e^{ \beta(\lambda) (\tau + L)} s_1(\tau + L, \lambda)
+ c_2 e^{-\beta(\lambda) (\tau + L)}
s_2( \tau + L, \lambda) = 0.
$$
Due to
$$
s_1(\tau, \lambda) = (-1)^N s_1(\tau + L, \lambda),~~~
 s_2(\tau, \lambda) = (-1)^N s_2(\tau + L, \lambda),~~~
$$
we must have
$$
c_1 s_1(\tau, \lambda) = 0, ~~ and ~~ c_2 s_2( \tau, \lambda) = 0.
$$
By the Sturm Separation Theorem\cite{eas2}, the zeros of
$s_1(x, \lambda) $ are separated from the zeros of $ s_2(x, \lambda) $.
Because  $s_1(\tau, \lambda) $ and $ s_2( \tau, \lambda) $ can not be zero
simultaneously, either $c_1$ or $c_2$ must be zero.
Whether we have a solution or not in this case is dependent on whether
either one of
$$
s_1(\tau, \lambda) = 0~~~~~~~~~~~~ and ~~~~~~~~~~~c_2 = 0
\eqno(E.1)
$$
and
$$
s_2( \tau, \lambda) = 0~~~~~~~~~~~~ and ~~~~~~~~~~~c_1 = 0
\eqno(E.2)
$$
is true.
From the discussion on (D1), (D2), (E), it can be seen that in the
conditional instability intervals in which $D(\lambda) \leq -2$, if we have
a solution of (2) and (3), we must have either $\mu_{2m} = \mu_{2m+1}$ or
one of (D2.1), (E.1), (E.2) when $\mu_{2m} < \mu_{2m+1}$.
From these equations we can see that we also always have
$y(\tau~+~a, \Lambda) = 0$ if $y(\tau, \Lambda) =0$.
Thus for $D(\lambda) \leq -2$, the same Equation (6) is a necessary
condition of (4).\par
It is easy to see that (6) is also a sufficient condition for having
a solution of Equations (2) and (3) for $|D(\lambda)| \geq 2$:
Equation (3) can be obtained by repeating (6) for $N$ times.\\
Thus {\it the Equation (6) is a necessary and sufficient condition
for having a solution of Equations (2) and (3) in a conditional instability
interval for which $ |D(\lambda)| \geq 2$}.
\par
Equation (6) does not contain the confinement length $L$, thus the eigen value
$\Lambda$ in a conditional instability interval
$ |D(\lambda)| \geq 2$ might be dependent on $\tau$, but not on $L$.
\par
The Theorem 3.1.3 of \cite{eas} indicates that
{\it for any real number $\tau$, there is always one and only one $\Lambda$
in each conditional
instability interval $[\mu_{2m}, \mu_{2m+1}]$ or $[\nu_{2m+1}, \nu_{2m+2}]$,
for which we have $y(\tau, \Lambda) = y(\tau~+~a, \Lambda) = 0$}.
We can label these $\Lambda$ as $\Lambda_{\tau, 2m}$ and
$\Lambda_{\tau, 2m+1}$.
According to the theorem 3.1.3 of e\cite{eas}, the range of
$\Lambda_{\tau, 2m}$ are in $[\mu_{2m}, \mu_{2m+1}]$ and
$\Lambda_{\tau, 2m+1}$ are in $[\nu_{2m+1}, \nu_{2m+2}]$.
\par
It is interesting to see how each of these $\tau$-dependent eigenvalues
$\Lambda$ changes as $\tau$ changes.  Here we only discuss
$\Lambda_{\tau, 2m}$. $\Lambda_{\tau, 2m+1}$ can be very similarly
discussed.
\par
If $\mu_{2m}~=~\mu_{2m+1}$, then we have the case D1,
the $\Lambda_{\tau, 2m}$ will not change as $\tau$ changes, we always
have $\Lambda_{\tau, 2m}~=~\mu_{2m}$ and the corresponding eigen function
$y(x, \Lambda_{\tau, 2m}) $ will be dependent on $\tau$, but it will always be
a semi-periodic function.
\par
The more interesting case is when $\mu_{2m}~\neq~\mu_{2m+1}$.
According to the theorem 3.1.2 of \cite{eas},
$\xi_{2m}(x)$ and $\xi_{2m+1}(x)$ have exact
$2m+1$ zeros in $[0,a)$.
Then according to the Sturm Comparison Theorem\cite{eas2}, the zeros of
$\xi_{2m}(x)$ and the zeros of $\xi_{2m+1}(x)$ must be
distributed alternatively:
There is always one and only one zero of $\xi_{2m+1}(x)$ between
two consecutive zeros of $\xi_{2m}(x)$, and
there is always one and only one zero of $\xi_{2m}(x)$ between
two consecutive zeros of $\xi_{2m+1}(x)$.
\par
We start from a zero $x_{1,2m}$ of $\xi_{2m}(x)$. If
$\tau~=~x_{1,2m}$, we have a solution of (2) and
(3) $\Lambda_{\tau,2m} = \mu_{2m}$,
and $y(x, \Lambda_{\tau,2m}) = \xi_{2m}(x)$.
There is a nearby zero $x_{1,2m+1}$ of $\xi_{2m+1}(x)$. If
$\tau = x_{1,2m+1}$, we have a solution of (2)
and (3) $\Lambda_{\tau,2m} = \mu_{2m+1}$,
and $y(x, \Lambda_{\tau,2m}) = \xi_{2m+1}(x)$.
If we treat $\tau$ as a variable and let $\tau$ go continuously from a zero
$x_{1,2m}$ of $\xi_{2m}(x)$ to a nearest neighbor zero $x_{1,2m+1}$
of $\xi_{2m+1}(x)$, because $\Lambda_{\tau, 2m}$ is a continuous
function of $\tau$\cite{eas}, the corresponding solution
of (2) and (3) $\Lambda_{\tau, 2m}$ will also go
continuously from $\mu_{2m}$ to $\mu_{2m+1}$.   Similarly,
if $\tau$ goes from the zero $x_{1,2m+1}$
of $\xi_{2m+1}(x)$ to the other nearest neighbor zero $x_{2,2m}$
of $\xi_{2m}(x)$ continuously, the corresponding $\Lambda_{\tau, 2m}$
will also go continuously from $\mu_{2m+1}$ to $\mu_{2m}$.
Since in the interval [0, a), both $\xi_{2m}(x)$ and $\xi_{2m+1}(x)$
have exactly $2m+1$ zeros, in general $\Lambda_{\tau, 2m}$ as function of
$\tau$ will always complete $2m+1$ ups and downs in an interval of length $a$.
Similarly $\Lambda_{\tau, 2m+1}$ as function of $\tau$ will
always complete $2m+2$ ups and downs in an interval of length $a$.
\par
Depending on the value of $\tau$, the solutions of (2) and (3) may
have different forms.
If $\tau$ goes in one direction (either in the right or in the left) from a
zero $x_{1,2m}$ of $\xi_{2m}(x)$ ($ y(x,\Lambda) =\xi_{2m}(x)$ here)
to its nearest neighbor zero $x_{1,2m+1}$ of $\xi_{2m+1}(x)$
($ y(x,\Lambda) =\xi_{2m+1}(x)$ now) then again to the next zero
$x_{2,2m}$ of $\xi_{2m}(x)$ ($ y(x,\Lambda) =\xi_{2m}(x)$ again)
continuously, the corresponding
$\Lambda_{\tau, 2m}$ will also goes from  $\mu_{2m}$ to $\mu_{2m+1}$
then again back to $\mu_{2m}$ continuously. But for $\tau$ in the two
open intervals ($x_{1,2m}, x_{1,2m+1}$) and ($x_{1,2m+1}, x_{2,2m}$) of this
path, the solution of (2) and (3) - the function $y(x, \Lambda)$ (4) has
two different forms: one has the form
$ c_1 e^{\beta(\Lambda)~x}~s_1(x, \Lambda)$,
and the other has the form $ c_2 e^{- \beta(\Lambda)~x}~s_2(x, \Lambda)$.
In which one section it has which form is dependent on $p(x),~q(x)$ and $s(x)$
in (1). \par
A function with the form of $ c_1 e^{\beta(\Lambda)~x}~s_1(x, \Lambda)$
or $ c_2 e^{- \beta(\Lambda)~x}~s_2(x, \Lambda)$, in which
$\beta(\Lambda) > 0$, is mainly distributed near either one of the two
ends of the truncated region, due to the exponential factor.
We can call it as an end solution of (2) and (3).
End solutions may also exist in the instability intervals $D (\Lambda) > 2$
($\nu_{2m+1} < \Lambda < \nu_{2m+2}$) and have
the form of either $ c_1 e^{\beta(\Lambda)~x}~p_1(x, \Lambda)$
or $ c_2 e^{- \beta(\Lambda)~x}~p_2(x, \Lambda)$ in the truncated region,
where $p_1(x, \Lambda)$ and $p_2(x, \Lambda)$ are period functions.
{\it These end solutions are introduced into the instability intervals
when the truncated boundary $\tau$ is not a zero of either one of the two
eigen functions $\zeta_n(x)$
of the periodic eigenvalue problem or when $\tau$ is not a zero of
either one of the two eigenfunctions $\xi_n(x)$
of the semi-periodic eigenvalue problem}. That
is the case $| D (\Lambda)| > 2$. \par
In Figure 1 is shown $\Lambda_{\tau,0}$ as the function of $\tau$.
In the open interval $(x_{1,0},x_{1,1})$ and
the open interval $(x_{1,1},x_{1,0}+a)$
the function $y(x, \Lambda)$ (4) has different
forms: one has the form $ c_1 e^{\beta(\Lambda)~x}~s_1(x, \Lambda)$,
and the other has the form $ c_2 e^{- \beta(\Lambda)~x}~s_2(x, \Lambda)$,
shown as a dashed line and a dotted line, indicating that the end
solution is mainly distributed near either one of the two different ends
of the truncated region. \par
In summary thus there are two different types of solutions of
the truncated Hill's equation: there are always
$N-1$ solutions in each stability interval, whose eigen values depend
on the truncation length $L$ but not on the truncation boundary $\tau$.
There is always one and only one solution in each finite conditional
instability interval $[\mu_{2m}, \mu_{2m+1}]$
or $[\nu_{2m+1},\nu_{2m+2}]$, whose eigen value might
depend on the boundary $\tau$ but not on the truncation length $L$.\par
\section{A special example: Electronic states in one dimensional crystals
of finite length}
\par
The one-dimensional Schr$\ddot{\rm o}$dinger differential equation\cite{kit}
with periodic potential is a special case of Hill's equation (1),
in which
\begin{equation}
p(x) = \frac{\hbar^2}{2m}, ~~s(x) = 1, ~~ q(x) = v(x)
\end{equation}
is the periodic potential and $\lambda = E $ is the energy eigen
value. Here $\hbar$ is the Planck's constant and $m$ is the
electron mass. It has been well known that the corresponding
eigenvalues form energy bands and the eigenfunctions are Bloch
waves\cite{kit}. The more general results obtained here for the
truncated Hill's equation can be directly applied to the
electronic states in one dimensional crystals of finite
length.
\par
By using a Kronig-Penney model potential, Pedersen and
Hemmer found that for a finite crystal of length $L = Na$, there
are $N-1$ states corresponding to each energy band $ n
> 0$, the energies of these $N-1$ states map the energy
bands\cite{ped} of the infinite crystal exactly and this mapping
does not depend on the boundaries of the finite crystal. This is a
special case of the more general result obtained in (5).\par By
using the theorem 3.1.2 of \cite{eas}, the author has given an
analytical solution on the complete quantum confinement of one
dimensional Bloch waves\cite{syr} in an {\it inversion-symmetric~
potential}. The results obtained in \cite{syr} greatly generalized
and improved many results obtained in \cite{ped}. \par The results
obtained in this work can be directly applied to the electronic
states in one dimensional crystals of finite length. The exact and
general
results are\cite{syr2}:\\
For one dimensional crystals bounded at $\tau$ and $\tau + L$,
there are two different types of electronic states: There are $N-1$ states
corresponding to each energy band of the Bloch wave,
their eigenvalues $\Lambda$ are given by (5), thus are dependent
only on the crystal length $L$
but not on the crystal boundary $\tau$ and map the energy band exactly;
There is always one and only one electronic state corresponding to
each band gap of the Bloch wave,
whose eigenvalue $\Lambda$ might be dependent on the crystal boundary
$\tau$ but not on the crystal length $L$.
Such a $L$-independent state can be either a constant energy
confined band edge state (if $\tau$ is a zero of a band edge wavefunction)
or a surface state in the band gap (if $\tau$ is not a zero of either band
edge wavefunction).
The results obtained in \cite{syr} are a special case of the more general
results obtained here.
\vskip .2in
Acknowledgment:
\par
The author is grateful to Professor Lo Yang for
stimulating discussions and Professor Shuxiang Yu for critical reading
of the manuscript. This research is supported by the National Natural
Science Foundation of China (Project No. 19774001).

\newpage
\begin{center}
Figure Captions
\end{center}

Fig. 1. $\Lambda_{\tau, 0}$ as a function of $\tau$ when $\mu_0~\neq~\mu_1$
in interval $x_{1,0} \leq \tau \leq x_{1,0} + a$.
The zeros of $\xi_{0}(x)$ are shown as solid circles and the
zero of $\xi_{1}(x)$ are shown as an open circle. \par
However, in general, $\Lambda_{\tau, n}$ as a function of $\tau$ does not
have to be linear.
   \\

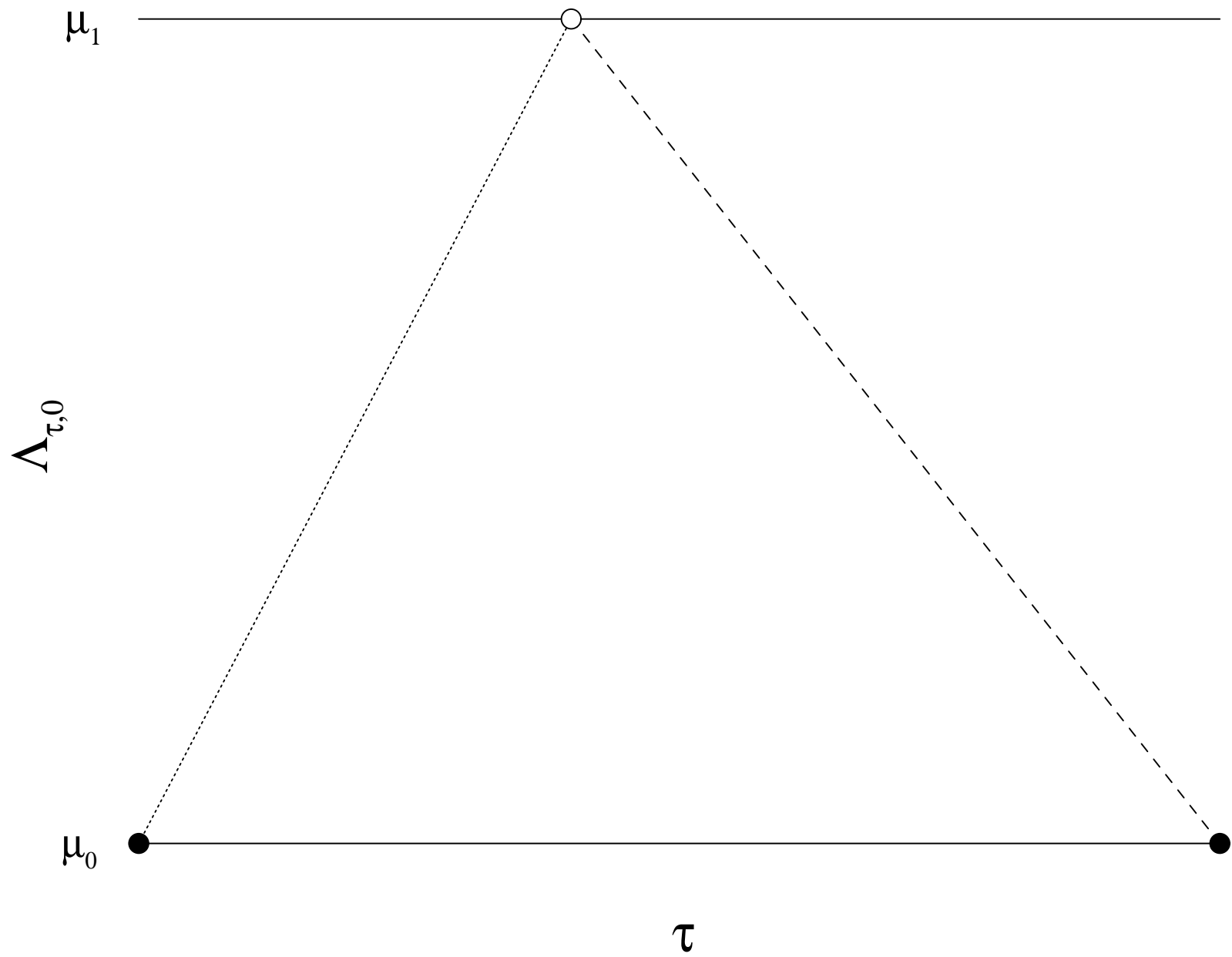
\end{document}